\documentclass[final,3p,times]{elsarticle}
\usepackage{hyperref}
\usepackage{amsmath}
\usepackage{amssymb}
\usepackage{amsthm}
\usepackage{listings}
\journal{Journal of Computational and Applied Mathematics}
\newtheorem{thm}{Theorem}
\newproof{prf}{Proof}
\begin{document}
\lstloadlanguages{C}
\lstset{language=C,extendedchars=false,numbers=left,numberstyle=\tiny,frame=lines,basicstyle=\small\ttfamily,columns=fullflexible,texcl=false,mathescape=true,lineskip=.25\baselineskip}
\begin{frontmatter}
  \title{Accurate complex Jacobi rotations\tnoteref{tit1}}
  \tnotetext[tit1]{The supplementary material is available in
    \href{https://github.com/venovako/AccJac}{\texttt{AccJac}} and
    \href{https://github.com/venovako/libpvn}{\texttt{libpvn}}
    repositories from \url{https://github.com/venovako} or by
    request.}
  \author[aut1]{Vedran Novakovi\'{c}\corref{cor1}}
  \ead{venovako@venovako.eu}
  \cortext[cor1]{corresponding author; \url{https://orcid.org/0000-0003-2964-9674}}
  \affiliation[aut1]{
    organization={independent researcher},
    addressline={Vankina ulica 15},
    city={Zagreb},
    postcode={HR-10020},
    country={Croatia}
  }
  \begin{abstract}\looseness=-1
    This note shows how to compute, to high relative accuracy under
    mild assumptions, complex Jacobi rotations for diagonalization of
    Hermitian matrices of order two, using the correctly rounded
    functions $\mathtt{cr\_hypot}$ and $\mathtt{cr\_rsqrt}$, proposed
    for standardization in the C programming language as recommended
    by the IEEE-754 floating-point standard.  The rounding to nearest
    (ties to even) and the non-stop arithmetic are assumed.  The
    numerical examples compare the observed with theoretical bounds on
    the relative errors in the rotations' elements, and show that the
    maximal observed departure of the rotations' determinants from
    unity is smaller than that of the transformations computed by
    LAPACK\@.
  \end{abstract}
\end{frontmatter}
\section{Introduction}\label{s:1}
Given $A$, a Hermitian matrix of order two, a complex Jacobi rotation
is a unitary matrix $U$ such that $\det{U}=1$ and $AU=U\Lambda$, i.e.,
$U^{\ast}AU=\Lambda$.  Here, $\Lambda$ is a real diagonal matrix of
the eigenvalues of $A$, and $U$ is a matrix of the associated
eigenvectors.  Let the representation of $A$, $U$, and $\Lambda$ be
fixed, with $\varphi\in[-\pi/4,\pi/4]$ and $\alpha\in(-\pi,\pi]$, as
\begin{equation}
  A=\begin{bmatrix}
  a_{11} & \overline{a_{21}}\\
  a_{21} & a_{22}
  \end{bmatrix},\quad
  U=\begin{bmatrix}
  \cos\varphi & -\mathrm{e}^{-\mathrm{i}\alpha}\sin\varphi\\
  \mathrm{e}^{\mathrm{i}\alpha}\sin\varphi & \cos\varphi
  \end{bmatrix},\quad
  \Lambda=\begin{bmatrix}
  \lambda_1 & 0\\
  0 & \lambda_2
  \end{bmatrix}.
  \label{e:1}
\end{equation}
It is shown here how to compute $\widetilde{U}$, a floating-point
approximation of $U$, from the rounded, floating-point value of $A$,
$\mathop{\mathrm{fl}}(A)$.  Empirically, such $\widetilde{U}$ is
closer, in the worst case, to a unitary matrix than the worst-case
transformation from the \texttt{xLAEV2}
LAPACK~\cite{Anderson-et-al-99} routine, as long as
$\mathop{\mathrm{fl}}(A)$ can be exactly scaled to $A'$, which has no
non-zero subnormal values.

Applications of the eigendecomposition (EVD) of a Hermitian matrix of
order two are numerous (see~\cite{Novakovic-23} for a more exhaustive
list), the most prominent being the Jacobi method for the EVD of
general Hermitian matrices~\cite{Hari-BegovicKovac-21}, and the
implicit (i.e., one-sided) Jacobi (also known as Hestenes) method for
the singular value decomposition (SVD) of $m\times n$, $m\ge n$
matrices, which is implemented in
LAPACK following~\cite{Drmac-97,Drmac-Veselic-08a,Drmac-Veselic-08b}
as \texttt{xGESVJ} and \texttt{xGEJSV} routines.

This work simplifies, and improves accuracy of, the complex Jacobi
rotations' floating-point computation from~\cite{Novakovic-23}, by
relying on \emph{correctly rounded} functions
$\mathop{\mathtt{cr\_hypot}}(x,y)=\mathop{\mathrm{fl}}(\!\sqrt{x^2+y^2})$
and
$\mathop{\mathtt{cr\_rsqrt}}(x)=\mathop{\mathrm{fl}}(1/\!\sqrt{x})$,
expected to be standardized as optional in C~\cite[\S7.33.8]{C-23}.
These functions have been provided\footnote{The source code and other
resources are freely available at
\url{https://core-math.gitlabpages.inria.fr} at the time of writing.}
by the CORE-MATH project~\cite{Sibidanov-et-al-22} for the default
rounding to nearest-even, in single (binary32) and double (binary64)
precision.  The latter function has also been considered
in~\cite{Borges-et-al-22}, and one approach towards the former has
been described in~\cite{Borges-20}.  It is therefore assumed here that
the non-stop floating-point arithmetic~\cite{IEEE-754-2019} is used,
with gradual underflow and rounding to nearest-even.

Usage of $\mathtt{hypot}$ for improving accuracy of real Jacobi
rotations has been considered in~\cite{Borges-17}, in a different
way and without the assumption on correct rounding.  As shown
in~\cite{Novakovic-23}, $\mathtt{cr\_rsqrt}$ can be safely avoided
when it is unavailable.

\looseness=-1
The special case of symmetric matrices and real Jacobi rotations is
not detailed for brevity, but the corresponding theoretical and
numerical results as well as the implementations can be be derived
from the presented ones, and are included in the supplement.
Performance of the code is not considered due to novelty of the
correctly rounded routines.
\section{Floating-point computation of complex Jacobi rotations}\label{s:2}
Listing~\ref{l:1} presents \texttt{ZJAEV2}\footnote{See
\url{https://github.com/venovako/libpvn/blob/master/src/pvn_ev2.c},
along with \texttt{pvn.h} and \texttt{pvn\_crm.h} in the same
subdirectory, and
\url{https://github.com/venovako/AccJac/blob/master/src/zjaev2.f90},
for the mostly LAPACK-compatible implementation.}, the proposed
algorithm for the Hermitian EVD in the form of~\eqref{e:1}, as a
simplification of~\cite[Listing~1]{Novakovic-23}.  In single
precision, \texttt{CJAEV2} can be obtained straightforwardly.  The
derivation of the algorithm can be found in~\cite{Novakovic-23},
including the scaling of $\mathop{\mathrm{fl}}(A)$ by a suitable power
of two to avoid overflow of the scaled eigenvalues.

Two accuracy-improving (but possibly performance-degrading)
modifications of~\cite{Novakovic-23} are proposed.  With $\nu$ being
the largest finite floating-point value, and using the C language
functions as defined in~\cite{C-23}, \cite[Eq.~(2.6)]{Novakovic-23}
has become
\begin{equation}
  \tan(2\tilde{\varphi})=\mathop{\text{\texttt{copysign}}}(\mathop{\text{\texttt{fmin}}}(\mathop{\text{\texttt{fmax}}}(\mathop{\mathrm{fl}}((2|\tilde{a}_{21}'|)/|\tilde{a}|),0),\nu),\tilde{a}),\quad
  \tilde{a}=\mathop{\mathrm{fl}}(a_{11}'-a_{22}'),\quad
  A'=\mathop{\mathrm{fl}}(2^{\zeta}\mathop{\mathrm{fl}}(A)),
  \label{e:2}
\end{equation}
with all roundings explicitly stated, and
$|\tilde{a}_{21}'|\mskip-2mu=\mskip-2mu\mathop{\mathtt{cr\_hypot}}(\Re{a_{21}'},\Im{a_{21}'})$.
Also, \cite[Eq.~(2.5)]{Novakovic-23} has been simplified as\footnote{If $\Lambda$ is not needed, $\sec{\tilde{\varphi}}=\mathop{\mathtt{cr\_hypot}}(\tan{\tilde{\varphi}},1)$, $\cos{\tilde{\varphi}}=1/\sec{\tilde{\varphi}}$, and $\sin{\tilde{\varphi}}=\tan{\tilde{\varphi}}/\sec{\tilde{\varphi}}$ could be used instead for possibly better accuracy.}
\begin{equation}
  \tan{\tilde{\varphi}}=\mathop{\mathrm{fl}}(\tan(2\tilde{\varphi})/\mathop{\mathrm{fl}}(1+\mathop{\mathtt{cr\_hypot}}(\tan(2\tilde{\varphi}),1))),\quad
  \sec^2{\tilde{\varphi}}=\mathop{\text{\texttt{fma}}}(\tan{\tilde{\varphi}},\tan{\tilde{\varphi}},1),\quad
  \cos{\tilde{\varphi}}=\mathop{\mathtt{cr\_rsqrt}}(\sec^2{\tilde{\varphi}}),
  \label{e:3}
\end{equation}
where $\mathop{\text{\texttt{fma}}}(x,y,z)=\mathop{\mathrm{fl}}(x\cdot y+z)$.
The symbols having a tilde indicate the computed approximations of the
exact values; e.g., $\tan{\tilde{\varphi}}$ approximates
$\tan\varphi$.  The elements of $\mathop{\mathrm{fl}}(A)$ are
considered exact for the purposes of this paper, i.e.,
$A=\mathop{\mathrm{fl}}(A)$.  For all but extremely badly scaled
matrices $A$ no element underflows when scaled as in~\eqref{e:2}, so
$A'$ is exact in such cases.

\begin{lstlisting}[float=hbtp,caption=\texttt{ZJAEV2}: the EVD of a double precision Hermitian matrix of order two in pseudo-C when rounding to nearest,label=l:1]
// input: $\!\mathop{\mathrm{fl}}(A)$ as $\!a_{11}^{}$, $\!a_{22}^{}$, $\!\Re{a_{21}^{}}$, $\!\Im{a_{21}^{}}$; output: $\!\cos{\tilde{\varphi}}$, $\!\cos{\tilde{\alpha}}\sin{\tilde{\varphi}}$, $\!\sin{\tilde{\alpha}}\sin{\tilde{\varphi}}$; $\!\tilde{\lambda}_1'$, $\!\tilde{\lambda}_2'$, $\!\zeta'$ or $\!\tilde{\lambda}_1^{}$, $\!\tilde{\lambda}_2^{}$
int $\zeta_{11}^{}$, $\zeta_{22}^{}$, $\zeta_{21}^{\Re}$, $\zeta_{21}^{\Im}$, $\zeta$, $\zeta'$, $\eta'$ = DBL_MAX_EXP - 3; double $\check{\mu}$ = DBL_TRUE_MIN, $\nu$ = DBL_MAX;
// determine $\zeta$ assuming all inputs are finite; avoid taking the exponent of 0
frexp(fmax(fabs($a_{11}$), $\check{\mu}$), &$\zeta_{11}$); frexp(fmax(fabs($a_{22}$), $\check{\mu}$), &$\zeta_{22}$);
frexp(fmax(fabs($\Re{a_{21}^{}}$), $\check{\mu}$), &$\zeta_{21}^{\Re}$); frexp(fmax(fabs($\Im{a_{21}^{}}$), $\check{\mu}$), &$\zeta_{21}^{\Im}$);
$\zeta$ = $\eta'$ - $\max\{\zeta_{11}^{},\zeta_{22}^{},\zeta_{21}^{\Re},\zeta_{21}^{\Im}\}$; $\zeta'$ = -$\zeta$; // $\max\{\cdots\}$ has to be expanded using the >= operator
// scale the input matrix $A$ by $2^{\zeta}$ as $A'=2^{\zeta}A$; the results might be inexact only due to underflow
$a_{11}'$ = scalbn($a_{11}^{}$, $\zeta$); $a_{22}'$ = scalbn($a_{22}^{}$, $\zeta$); $\Re{a_{21}'}$ = scalbn($\Re{a_{21}^{}}$, $\zeta$); $\Im{a_{21}'}$ = scalbn($\Im{a_{21}^{}}$, $\zeta$);
// find the polar form of $a_{21}'=2^{\zeta}a_{21}^{}$ as $|\tilde{a}_{21}'|\mathrm{e}^{\mathrm{i}\tilde{\alpha}}=|\tilde{a}_{21}'|(\cos{\tilde{\alpha}}+\mathrm{i}\cdot\sin{\tilde{\alpha}})$; $|\Re{\tilde{a}_{21}'}|=|\Im{\tilde{a}_{21}'}|=\check{\mu}$ is dangerous
$|\Re{a_{21}'}|$ = fabs($\Re{a_{21}'}$); $|\Im{a_{21}'}|$ = fabs($\Im{a_{21}'}$); $|\tilde{a}_{21}'|$ = $\mathtt{cr\_hypot}$($|\Re{a_{21}'}|$, $|\Im{a_{21}'}|$);
$\cos{\tilde{\alpha}}$ = copysign(fmin($|\Re{a_{21}'}|$ / $|\tilde{a}_{21}'|$, 1.0), $\Re{a_{21}'}$); $\sin{\tilde{\alpha}}$ = $\Im{a_{21}'}$ / fmax($|\tilde{a}_{21}'|$, $\check{\mu}$); // underflows possible
// compute the complex Jacobi rotation
$\tilde{o}$ = $|\tilde{a}_{21}'|$ * 2.0; $\tilde{a}$ = $a_{11}'$ - $a_{22}'$; $|\tilde{a}|$ = fabs($a$); // underflow of $\tilde{a}$ possible
$\tan(2\tilde{\varphi})$ = copysign(fmin(fmax($\tilde{o}$ / $|\tilde{a}|$, 0.0), $\nu$), $\tilde{a}$); // filter out possible non-finite partial results
$\tan{\tilde{\varphi}}$ = $\tan(2\tilde{\varphi})$ / (1.0 + $\mathtt{cr\_hypot}$($\tan(2\tilde{\varphi})$, 1.0)); // rounding to nearest: denominator cannot overflow
$\sec^2{\tilde{\varphi}}$ = fma($\tan{\tilde{\varphi}}$, $\tan{\tilde{\varphi}}$, 1.0); $\cos{\tilde{\varphi}}$ = $\mathtt{cr\_rsqrt}$($\sec^2{\tilde{\varphi}}$); $\sin{\tilde{\varphi}}$ = $\tan{\tilde{\varphi}}$ * $\cos{\tilde{\varphi}} $;
$\cos{\tilde{\alpha}}\sin{\tilde{\varphi}}$ = $\cos{\tilde{\alpha}}$ * $\sin{\tilde{\varphi}}$; $\sin{\tilde{\alpha}}\sin{\tilde{\varphi}}$ = $\sin{\tilde{\alpha}}$ * $\sin{\tilde{\varphi}}$; // $\tan{\tilde{\varphi}}$, $\sin{\tilde{\varphi}}$, $\Re(\mathrm{e}^{\mathrm{i}\tilde{\alpha}}\sin{\tilde{\varphi}})$, $\Im(\mathrm{e}^{\mathrm{i}\tilde{\alpha}}\sin{\tilde{\varphi}})$ can underflow
// compute the eigenvalues scaled by $2^{\zeta}$ (always finite)
$\tilde{\lambda}_1'$ = fma($\tan{\tilde{\varphi}}$, fma($a_{22}'$, $\tan{\tilde{\varphi}}$, $\tilde{o}$), $a_{11}'$) / $\sec^2{\tilde{\varphi}}$; $\tilde{\lambda}_2'$ = fma($\tan{\tilde{\varphi}}$, fma($a_{11}'$, $\tan{\tilde{\varphi}}$, -$\tilde{o}$), $a_{22}'$) / $\sec^2{\tilde{\varphi}}$;
$\tilde{\lambda}_1^{}$ = scalbn($\tilde{\lambda}_1'$, $\zeta'$); $\tilde{\lambda}_2^{}$ = scalbn($\tilde{\lambda}_2'$, $\zeta'$); // optionally backscale (over/under-flows possible)
\end{lstlisting}

\texttt{ZJAEV2}, unlike \texttt{ZLAEV2}, does not attempt to sort the
computed eigenvalues decreasingly, leaving it to users instead.
\section{Rounding error analysis}\label{s:3}
\looseness=-1
Theorem~\ref{t:1} shows that the elements of
$\widetilde{U}$ computed by \texttt{xJAEV2} have high relative
accuracy, with the error bounds tighter than
in~\cite[Proposition~2.4]{Novakovic-23}, if no computed inexact value
underflows.  Exact subnormals and zeros are allowed.

\begin{thm}\label{t:1}
  In at least single precision arithmetic with rounding to nearest,
  and barring any inexact underflow, the elements of $\widetilde{U}$
  computed by \texttt{xJAEV2} have the relative error bounds as
  specified, with $\varepsilon$ being the machine precision:
  \begin{equation}
    \begin{gathered}
      \cos{\tilde{\varphi}}=\delta_c^{}\cos\varphi,\quad
      1-6.00000017\,\varepsilon<\delta_c^{}<1+6.00000000\,\varepsilon,\\
      \cos{\tilde{\alpha}}\sin{\tilde{\varphi}}=\delta_{21}^{\Re}\cos\alpha\sin\varphi,\ \ \
      \sin{\tilde{\alpha}}\sin{\tilde{\varphi}}=\delta_{21}^{\Im}\sin\alpha\sin\varphi,\quad
      1-19.00000000\,\varepsilon<\delta_{21}^{\Re},\delta_{21}^{\Im}<1+19.00000950\,\varepsilon.
    \end{gathered}
    \label{e:4}
  \end{equation}
\end{thm}
\begin{prf}\looseness=-1
  Since $\mathtt{cr\_hypot}$ is correctly rounded,
  $|\tilde{a}_{21}'|=|a_{21}'|(1+\epsilon_1^{})$,
  $|\epsilon_1^{}|\le\varepsilon$, and
  $|\tilde{a}_{21}'|=0$ if and only if $a_{21}'=0$.  With
  $\epsilon_2^{}$ and $\epsilon_3^{}$,
  $\max\{|\epsilon_2^{}|,|\epsilon_3^{}|\}\le\varepsilon$, from
  line~11 of Listing~\ref{l:1}, when $|\tilde{a}_{21}'|>0$, the
  relative errors $\delta_{\alpha}^{\Re}$ and $\delta_{\alpha}^{\Im}$
  in $\cos{\tilde{\alpha}}$ and $\sin{\tilde{\alpha}}$ are
  \begin{displaymath}
    \cos{\tilde{\alpha}}=\frac{\Re{a_{21}'}}{|a_{21}'|(1+\epsilon_1^{})}(1+\epsilon_2^{})=\cos\alpha\frac{1+\epsilon_2^{}}{1+\epsilon_1^{}}=\delta_{\alpha}^{\Re}\cos\alpha,\quad
    \sin{\tilde{\alpha}}=\frac{\Im{a_{21}'}}{|a_{21}'|(1+\epsilon_1^{})}(1+\epsilon_3^{})=\sin\alpha\frac{1+\epsilon_3^{}}{1+\epsilon_1^{}}=\delta_{\alpha}^{\Im}\sin\alpha.
  \end{displaymath}
  The minimal value of $(1+\epsilon_i^{})/(1+\epsilon_1^{})$,
  $i\in\{2,3\}$, is
  $\delta_{\alpha}^-=(1-\varepsilon)/(1+\varepsilon)$, and the maximal
  value is $\delta_{\alpha}^+=(1+\varepsilon)/(1-\varepsilon)$.
  When $\tilde{a}_{21}'$ is real or imaginary (let $\beta=1$ then, and
  otherwise $\beta=2$), $\cos{\tilde{\alpha}}$ and
  $\sin{\tilde{\alpha}}$ are exact, so
  $\delta_{\alpha}^-=\delta_{\alpha}^+=1$.  The scaling in lines~4--8
  of Listing~\ref{l:1} prevents $|\tilde{a}_{21}'|$, as well as
  $\tilde{o}$, from overflowing (see~\cite[Proposition~2.4 and
    Corollary~2.5]{Novakovic-23}).

  Multiplication by a power of two is exact if the result is normal,
  so $\tilde{o}=o(1+\epsilon_1)^{\beta-1}$ in line~13 of
  Listing~\ref{l:1}.  It is assumed that $\tilde{a}$ is exact (i.e.,
  $\tilde{a}=a$) or has not underflown, so
  $|\tilde{a}|=|a|(1+\epsilon_4)$, $|\epsilon_4|\le\varepsilon$.  With
  $|\epsilon_5|\le\varepsilon$ it holds
  \begin{equation}
    \mathop{\mathrm{fl}}\left(\frac{\tilde{o}}{|\tilde{a}|}\right)=\frac{o(1+\epsilon_1^{})^{\beta-1}}{|a|(1+\epsilon_4^{})}(1+\epsilon_5^{})=\delta_2^{}\frac{o}{|a|},\quad
    \delta_2^{}=\frac{(1+\epsilon_1^{})^{\beta-1}(1+\epsilon_5^{})}{1+\epsilon_4^{}},\quad
    \frac{(1-\varepsilon)^{\beta}}{1+\varepsilon}=\delta_2^-\le\delta_2^{}\le\frac{(1+\varepsilon)^{\beta}}{1-\varepsilon}=\delta_2^+.
    \label{e:5}
  \end{equation}
  If $\tilde{o}=0$, then $\tan(2\tilde{\varphi})=\pm 0$, so it is
  exact.  Else, if $\tilde{a}=\pm 0$, then
  $\tan{\tilde{\varphi}}=\pm 1$ is exact, regardless of the error in
  $\tan(2\tilde{\varphi})$.  Else, if $\tilde{a}$ causes overflow of
  the result of the division in line~14 of Listing~\ref{l:1},
  $\tan(2\tilde{\varphi})=\pm\nu$ instead of $\pm\infty$.  Its
  relative error is then not easily bounded, but is irrelevant, since
  $\tan\varphi$ is a monotonically non-decreasing function of
  $\tan(2\varphi)$, as well as $\tan{\tilde{\varphi}}$ with respect to
  $\tan(2\tilde{\varphi})$, due to monotonicity of all floating-point
  operations involved in computing it (addition, division, and
  hypotenuse, with the associated roundings) as in~\eqref{e:3}.  The
  relative error in $\tan{\tilde{\varphi}}$ largest by magnitude,
  i.e., $\tan{\tilde{\varphi}}=\pm 1$, computed when
  $|\tan(2\varphi)|\le\nu$ for which $|\tan(2\tilde{\varphi})|=\nu$,
  can only decrease for $\tan(2\varphi)$ larger in magnitude, since
  $\tan(2\varphi)\to\pm\infty$ implies $\tan\varphi\to\pm 1$.
  Moreover, $\tan(2\tilde{\varphi})$ is not an output value of
  \texttt{xJAEV2}.

  The equation
  $\tan^2(2\tilde{\varphi})+1=(\tan^2(2\varphi)+1)(1+\epsilon_6)$ has
  to be solved for $\epsilon_6$, to relate the exact with the computed
  hypotenuse in line~15 of Listing~\ref{l:1}.  By substituting
  $\delta_2\tan(2\varphi)$ for $\tan(2\tilde{\varphi})$, what can be
  done when $|\tan(2\varphi)|\le\nu$, and rearranging the terms, while
  noting that $0<\delta_2^-<1$, $\delta_2^+>1$, and $0\le x/(x+1)<1$
  for all $x\ge 0$, it follows that
  \begin{equation}
    \epsilon_6^{}=\frac{\tan^2(2\varphi)}{\tan^2(2\varphi)+1}(\delta_2^2-1),\quad
    (\delta_2^-)^2-1<\epsilon_6^{}<(\delta_2^+)^2-1,\quad
    \delta_2^-<\!\sqrt{1+\epsilon_6^{}}<\delta_2^+.
    \label{e:6}
  \end{equation}
  Let $z=\!\sqrt{\tan^2(2\varphi)+1}$.  Then,
  $\!\sqrt{\tan^2(2\tilde{\varphi})+1}=z\sqrt{1+\epsilon_6}$, so
  $\tilde{z}=\mathop{\mathtt{cr\_hypot}}(\tan(2\tilde{\varphi}),1)=z\sqrt{1+\epsilon_6}(1+\epsilon_7)=z\cdot y$,
  where $|\epsilon_7|\le\varepsilon$.  Now,
  $1+\tilde{z}=(1+z)(1+\epsilon_8)$ has to be solved for $\epsilon_8$,
  to relate the exact and the computed denominator in line~15 of
  Listing~\ref{l:1}.  Similarly as above, since $z\ge 0$ and, due
  to~\eqref{e:6},
  $\delta_2^-(1-\varepsilon)=y^-<y<y^+=\delta_2^+(1+\varepsilon)$, it
  follows
  \begin{equation}
    \epsilon_8=\frac{z}{z+1}(y-1),\quad
    y^--1<\epsilon_8<y^+-1,\quad
    y^-<1+\epsilon_8<y^+.
    \label{e:7}
  \end{equation}
  Therefore,
  $\mathop{\mathrm{fl}}(1+\tilde{z})=(1+z)(1+\epsilon_8)(1+\epsilon_9)$,
  where $|\epsilon_9|\le\varepsilon$.  Let $x=\tan(2\varphi)$ and
  $\tilde{x}=\tan(2\tilde{\varphi})$.  Then, with
  $|\epsilon_{10}|\le\varepsilon$,
  \begin{equation}
    \tan{\tilde{\varphi}}=\mathop{\mathrm{fl}}\left(\frac{\tilde{x}}{\mathop{\mathrm{fl}}(1+\tilde{z})}\right)=\frac{x}{1+z}\frac{\delta_2(1+\epsilon_{10})}{(1+\epsilon_8)(1+\epsilon_9)}=\delta_t\tan\varphi,\quad
    \delta_t=\frac{\delta_2}{1+\epsilon_8}\frac{1+\epsilon_{10}}{1+\epsilon_9}.
    \label{e:8}
  \end{equation}
  It is easy to minimize the second factor of $\delta_t^{}$
  in~\eqref{e:8} to $(1-\varepsilon)/(1+\varepsilon)$, and to maximize
  it to $(1+\varepsilon)/(1-\varepsilon)$, irrespectively of the first
  factor.  Using~\eqref{e:7}, $1+\epsilon_8^{}$ is minimized to
  $\delta_2^-(1-\varepsilon)$ and maximized to
  $\delta_2^+(1+\varepsilon)$.  Therefore, from~\eqref{e:5} follows
  \begin{equation}
    \delta_t^-=\frac{\delta_2^-}{\delta_2^+(1+\varepsilon)}\frac{1-\varepsilon}{1+\varepsilon}=\frac{(1-\varepsilon)^{\beta+2}}{(1+\varepsilon)^{\beta+3}},\quad
    \delta_t^+=\frac{\delta_2^+}{\delta_2^-(1-\varepsilon)}\frac{1+\varepsilon}{1-\varepsilon}=\frac{(1+\varepsilon)^{\beta+2}}{(1-\varepsilon)^{\beta+3}},\quad
    \delta_t^-<\delta_t^{}<\delta_t^+.
    \label{e:9}
  \end{equation}

  \looseness=-1
  Similarly as in~\eqref{e:6}, the equation
  $\delta_t^2\tan^2\varphi+1=(\tan^2\varphi+1)(1+\epsilon_{11}^{})$
  has to be solved for $\epsilon_{11}^{}$, what, due to~\eqref{e:9},
  leads to
  \begin{equation}
    \epsilon_{11}^{}=\frac{\tan^2\varphi}{\tan^2\varphi+1}(\delta_t^2-1),\quad
    0\le\frac{\tan^2\varphi}{\tan^2\varphi+1}\le\frac{1}{2},\quad
    \frac{1}{\sqrt{2}}\sqrt{(\delta_t^-)^2+1}=\delta_r^-<\!\sqrt{1+\epsilon_{11}^{}}<\delta_r^+=\frac{1}{\sqrt{2}}\sqrt{(\delta_t^+)^2+1}.
    \label{e:10}
  \end{equation}
  Then, from line~16 of Listing~\ref{l:1},
  $\sec^2{\tilde{\varphi}}=(\delta_t^2\tan^2\varphi+1)(1+\epsilon_{12}^{})=\delta_q^2\sec^2\varphi$,
  where $|\epsilon_{12}^{}|\le\varepsilon$ and
  $\delta_q^2=(1+\epsilon_{11}^{})(1+\epsilon_{12}^{})$,
  while, from~\eqref{e:10},
  $\delta_r^-\sqrt{1-\varepsilon}=\delta_q^-<\delta_q^{}<\delta_q^+=\delta_r^+\sqrt{1+\varepsilon}$.
  With $|\epsilon_{13}^{}|\le\varepsilon$, the relative error in
  $\cos{\tilde{\varphi}}$ can be bounded as
  \begin{equation}
    \cos{\tilde{\varphi}}=\frac{1+\epsilon_{13}^{}}{\sqrt{\sec^2{\tilde{\varphi}}}}=\frac{1+\epsilon_{13}^{}}{\delta_q^{}}\cos\varphi=\delta_c^{}\cos\varphi,\quad
    \frac{1-\varepsilon}{\delta_q^+}=\delta_c^-<\delta_c^{}<\delta_c^+=\frac{1+\varepsilon}{\delta_q^-}.
    \label{e:11}
  \end{equation}

  Now,
  $\sin{\tilde{\varphi}}=\delta_s^{}\sin\varphi$, where
  $\delta_s^{}=\delta_t^{}\delta_c^{}(1+\epsilon_{14}^{})$,
  $|\epsilon_{14}^{}|\le\varepsilon$, and
  $\delta_t^-\delta_c^-(1-\varepsilon)=\delta_s^-<\delta_s^{}<\delta_s^+=\delta_t^+\delta_c^+(1+\varepsilon)$.
  It holds
  \begin{equation}
    \begin{gathered}
      \cos{\tilde{\alpha}}\sin{\tilde{\varphi}}=\delta_{21}^{\Re}\cos\alpha\sin\varphi,\quad
      \delta_{21}^{\Re}=\delta_{\alpha}^{\Re}\delta_s^{}(1+\epsilon_{15}^{})^{\beta-1},\quad
      \sin{\tilde{\alpha}}\sin{\tilde{\varphi}}=\delta_{21}^{\Im}\sin\alpha\sin\varphi,\quad
      \delta_{21}^{\Im}=\delta_{\alpha}^{\Im}\delta_s^{}(1+\epsilon_{16}^{})^{\beta-1},\\
      \max\{|\epsilon_{15}^{}|,|\epsilon_{16}^{}|\}\le\varepsilon,\quad
      \delta_{\alpha}^-\delta_s^-(1-\varepsilon)^{\beta-1}=\delta_{21}^-<\delta_{21}^{\Re},\delta_{21}^{\Im}<\delta_{21}^+=\delta_{\alpha}^+\delta_s^+(1+\varepsilon)^{\beta-1}.
    \end{gathered}
    \label{e:12}
  \end{equation}

  Given $p$, the number of non-implied bits of the significand of a
  floating-point value ($p=23,52,112$ for the binary32, binary64, and
  binary128 types, respectively), and thus $\varepsilon=2^{-p-1}$, the
  bounds $\delta_c^-$ and $\delta_c^+$ from~\eqref{e:11}, and
  $\delta_{21}^-$ and $\delta_{21}^+$ from~\eqref{e:12} were evaluated
  symbolically by the Wolfram Language
  script\footnote{\url{https://github.com/venovako/AccJac/blob/master/etc/thm1.wls},
  with the command-line arguments $p$, $\mathtt{n}$, $\beta$, and
  $\mathtt{d}$} from the supplement.  The expressions $1-\delta_c^-$,
  $\delta_c^+-1$, $1-\delta_{21}^-$, and $\delta_{21}^+-1$ were then
  divided by $\varepsilon$, approximated to $\mathtt{n}=100$ digits,
  and rounded away from zero to $\mathtt{d}=8$ decimal places, to
  calculate by how many $\varepsilon$ each bound differs from unity.
  This was repeated for all three aforementioned datatypes, although
  only for the first two the required correctly rounded functions
  exist at present, and could have been similarly done for another
  (e.g., half) precision.  Among the three bounds of the same (e.g.,
  $\delta_c^-$) kind, the one with the largest respective multiple of
  $\varepsilon$ became the overall worst-case bound of its kind (the
  lower bound for $\delta_c^{}$ in this case).  The worst-case bounds
  were finally collected in~\eqref{e:4}, what concludes the proof for
  the standard datatypes.
  \qed
\end{prf}
Even tighter error bounds hold when $a_{21}'$ is real or imaginary, as
indicated in the proof of Theorem~\ref{t:1} (take $\beta=1$), and can
as well hold for one or more datatypes in the general ($\beta=2$)
case, as shown in~\eqref{e:13} for double precision ($p=52$),
\begin{equation}
  \begin{array}{cccc}
    \delta_c^-\gtrapprox  1- 6.00000001\varepsilon,&
    \delta_c^+\lessapprox 1+ 6.00000000\varepsilon,&
    \delta_{21}^-\gtrapprox  1-19.00000000\varepsilon,&
    \delta_{21}^+\lessapprox 1+19.00000001\varepsilon.
  \end{array}
  \label{e:13}
\end{equation}
\section{Numerical testing}\label{s:4}
\looseness=-1
Let, for an exact quantity $x\ne 0$, the relative error in an
approximation $\tilde{x}$ in the terms of the machine precision
$\varepsilon$ be
\begin{displaymath}
  \mathop{\rho}(\tilde{x})=(\tilde{x}-x)/(x\cdot\varepsilon).
\end{displaymath}
If $x=0$, then $\mathop{\rho}(\tilde{x})=0$ if $\tilde{x}=0$ and
$\pm\infty$ otherwise.  For a computed Jacobi rotation
$\widetilde{U}$, take as a measure of its non-unitarity the difference
of its determinant from unity, in the terms of $\varepsilon$ to keep
the measure's values comparable across the spectrum of floating-point
datatypes.  Since
$\tilde{u}_{11}=\tilde{u}_{22}=\cos{\tilde{\varphi}}$ and
$\tilde{u}_{21}=-\overline{\tilde{u}_{12}}=\cos{\tilde{\alpha}}\sin{\tilde{\varphi}}+\mathrm{i}\cdot\sin{\tilde{\alpha}}\sin{\tilde{\varphi}}$,
define
\begin{displaymath}
  \mathop{\Delta}(\widetilde{U})=\mathop{\rho}(\det\widetilde{U})=(\det\widetilde{U}-1)/\varepsilon=(\cos^2{\tilde{\varphi}}+(\cos{\tilde{\alpha}}\sin{\tilde{\varphi}})^2+(\sin{\tilde{\alpha}}\sin{\tilde{\varphi}})^2-1)/\varepsilon.
\end{displaymath}

Each of $33$ test runs generated random Hermitian matrices $A_k$,
$1\le k\le 2^{30}$, and computed
$\widetilde{U}_{k;\text{\texttt{Z}}}^{[\text{\texttt{X}}]}$,
$\Delta_{\text{\texttt{X}}}^-$, $\Delta_{\text{\texttt{X}}}^+$,
$\mathop{\rho_{\text{\texttt{J}}}^-}(\cos{\tilde{\varphi}})$,
$\mathop{\rho_{\text{\texttt{J}}}^+}(\cos{\tilde{\varphi}})$,
$\mathop{\rho_{\text{\texttt{J}}}^-}(\Re{\tilde{u}_{21}^{}})$,
$\mathop{\rho_{\text{\texttt{J}}}^+}(\Re{\tilde{u}_{21}^{}})$,
$\mathop{\rho_{\text{\texttt{J}}}^-}(\Im{\tilde{u}_{21}^{}})$, and
$\mathop{\rho_{\text{\texttt{J}}}^+}(\Im{\tilde{u}_{21}^{}})$, where
$\text{\texttt{X}}=\text{\texttt{J}}$ stands for \texttt{ZJAEV2} and
$\text{\texttt{X}}=\text{\texttt{L}}$ for \texttt{ZLAEV2}, and the
superscripts $-$ and $+$ denote the minimal and the maximal relative
error in a particular quantity in the run; i.e.,
\begin{equation}
  \Delta_{\text{\texttt{X}}}^-=\min_k\mathop{\Delta}(\widetilde{U}_{k;\text{\texttt{Z}}}^{[\text{\texttt{X}}]}),\quad
  \Delta_{\text{\texttt{X}}}^+=\max_k\mathop{\Delta}(\widetilde{U}_{k;\text{\texttt{Z}}}^{[\text{\texttt{X}}]}),\quad
  \mathop{\rho_{\text{\texttt{J}}}^-}(\tilde{x})=\min_k\mathop{\rho}(\tilde{x}_k^{}),\quad
  \mathop{\rho_{\text{\texttt{J}}}^+}(\tilde{x})=\max_k\mathop{\rho}(\tilde{x}_k^{}),\quad
  \tilde{x}\in\{\cos{\tilde{\varphi}},\Re{\tilde{u}_{21}^{}},\Im{\tilde{u}_{21}^{}}\}.
  \label{e:14}
\end{equation}

The error values were obtained in quadruple precision (the MPFR
library~\cite{Fousse-et-al-07} could have been used for higher
precision).  For $\rho_{\text{\texttt{J}}}^{\pm}$ from~\eqref{e:14},
$\widetilde{U}_{k;\text{\texttt{Q}}}^{[\text{\texttt{J}}]}$ was
computed by \texttt{QJAEV2}, using the \texttt{\_\_float128} C
datatype and the corresponding functions from the GCC's
\texttt{libquadmath}\footnote{Some of its functions might not be
correctly rounded, as noted in~\cite{Lefevre-et-al-23} for
\texttt{sqrtq}, but this was not checked with the present
implementation.} library (substituting \texttt{hypotq} for
$\mathtt{cr\_hypot}$ and \texttt{1/sqrtq} for $\mathtt{cr\_rsqrt}$,
and calculating $\sin\tilde{\varphi}$
as$\tan\tilde{\varphi}/\!\sqrt{\sec^2\tilde{\varphi}}$), and was taken
to be the exact $U_k^{}$ from~\eqref{e:1}, while the inexact was
$\widetilde{U}_{k;\text{\texttt{Z}}}^{[\text{\texttt{J}}]}$.

The testing was performed on an Intel Xeon Phi 7210 machine, with the
GNU C (\texttt{gcc}) and Fortran (\texttt{gfortran}) compilers,
version 13.2.1, on a 64-bit Linux.  Random numbers were harvested from
\texttt{/dev/random} and discarded if not within $[\mu,\nu/4]$ by
magnitude, with $\mu$ being the smallest positive normal value.  Each
Hermitian matrix was defined by four such random numbers, two for the
offdiagonal element and two for the diagonal.  The reference
LAPACK\footnote{See \url{https://github.com/Reference-LAPACK} and
\url{https://github.com/venovako/AccJac/blob/master/src/zlaev2.f90}.}
routines were compiled from the modernized but semantically faithful
sources, because it is not easy to determine exactly how the optimized
routines in vendors' libraries compute nor whether or when their
behavior will change.

\begin{figure}[h!bt]
  \centering\includegraphics{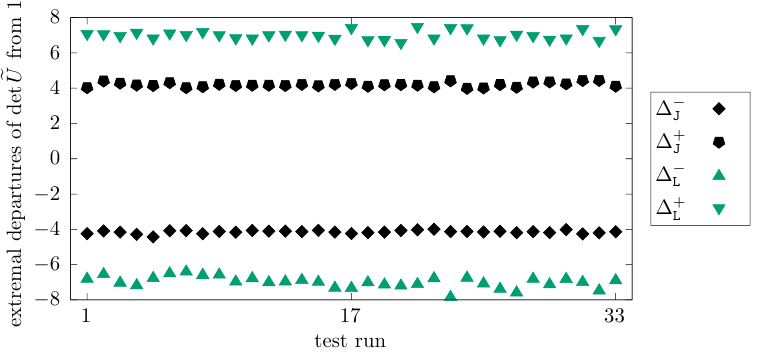}
  \caption{Extremal departures of $\det\widetilde{U}$ from unity in
    multiples of $\varepsilon$, for \texttt{ZJAEV2} and
    \texttt{ZLAEV2}.}
  \label{f:1}
\end{figure}

\looseness=-1
Figure~\ref{f:1} shows the values of
$\Delta_{\text{\texttt{X}}}^-$ and $\Delta_{\text{\texttt{X}}}^+$,
between which all
$\mathop{\Delta}(\widetilde{U}_{k;\text{\texttt{Z}}}^{[\text{\texttt{X}}]})$
in a particular run lie.  It can be concluded that in the worst case
\texttt{ZJAEV2} outputs a transformation that is almost twice
``closer'' to a unitary matrix than the worst-case one from
\texttt{ZLAEV2} as long as the elements of $A'$ are normal (otherwise,
see the comments on accuracy in~\cite[Sect.~5.2]{Novakovic-23}).

\begin{figure}[h!bt]
  \centering\includegraphics{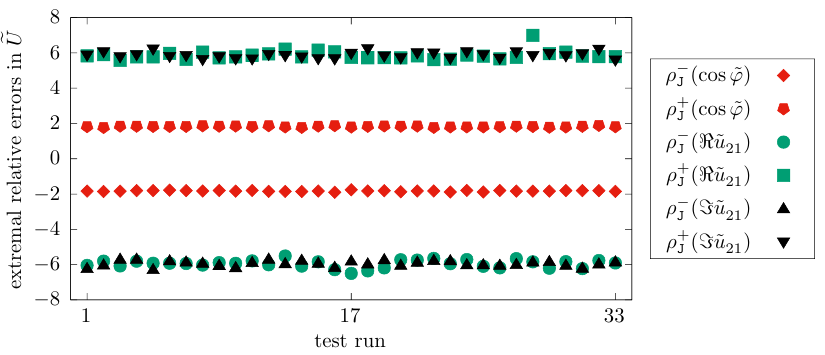}
  \caption{Observed extremal relative errors, in multiples of
    $\varepsilon$, in the elements of $\widetilde{U}$ computed by
    \texttt{ZJAEV2}.}
  \label{f:2}
\end{figure}

Figure~\ref{f:2} depicts the observed relative errors in the elements
of $\widetilde{U}$ computed by \texttt{ZJAEV2}, which fall well within
the ranges from~\eqref{e:4}.  This shows that the bounds from
Theorem~\ref{t:1} can probably be tightened, what is left for future
work.

\begin{figure}[h!bt]
  \centering\includegraphics{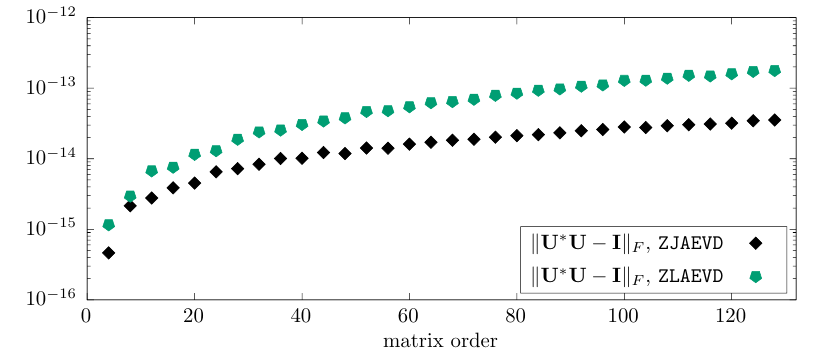}
  \caption{Departures from unitarity of the accumulated eigenvector
    matrices $\mathbf{U}$ of order $n$ from \texttt{ZJAEVD} and
    \texttt{ZLAEVD}.}
  \label{f:3}
\end{figure}

As a practical application of the proposed algorithm, a Jacobi
eigensolver for Hermitian matrices~\cite{Hari-BegovicKovac-21} of any
order $n$ has been implemented, that annihilates a pair of pivot
elements of the largest magnitude in each sequential step by
pre-multiplying the pivot rows of the iteration matrix by
$\widetilde{U}^{\ast}$ and post-multiplying the pivot columns of the
iteration matrix and of the eigenvector matrix $\mathbf{U}$ by
$\widetilde{U}$, where $\widetilde{U}$ is computed either by
\texttt{ZJAEV2} or by \texttt{ZLAEV2}. The former variant is called
\texttt{ZJAEVD} and the latter \texttt{ZLAEVD}.  Figure~\ref{f:3}
shows that the final $\mathbf{U}$ is closer to a unitary matrix with
\texttt{ZJAEV2} than with \texttt{ZLAEV2}, for the given Hermitian
matrices with eigenvalues $i$, $1\le i\le n$, where $n=4j$,
$1\le j\le 32$.

\appendix
\section*{Acknowledgements}\label{s:a}
Most of the computing resources used have remained available to the
author after the project IP--2014--09--3670 ``Matrix Factorizations
and Block Diagonalization Algorithms''
(\href{https://web.math.pmf.unizg.hr/mfbda}{MFBDA}) by Croatian
Science Foundation expired.

The author would like to thank Dean Singer for his material support
and declares no competing interests.
\end{document}